\input amstex
\documentstyle {amsppt}
\pageheight{50.5pc}
\pagewidth{32pc}
\define\mbR{\Bbb R}
\define\mfM{\frak M}
\define\ve{\varepsilon}
\define\vf{\varphi}
\define\vk{\varkappa}
\define\cF{\Cal F}
\define\cU{\Cal U}
\define\cG{\Cal G}

\topmatter
\title{}
{
ONE BROWNIAN STOCHASTIC FLOW
}
\endtitle
\author
Andrey A.Dorogovtsev
\endauthor
\address
Institute of Mathematics, Ukrainian Academy of
Science, ul. Tereshchenkivska, 3, 01601, Kiev, Ukraine
\endaddress
\keywords
Stochastic flows, measure-valued processes, weak compactness, Brownian
motion
\endkeywords
\subjclass
\endsubjclass
\email
adoro\@imath.kiev.ua
\endemail
\abstract
The weak limits of the measure-valued processes organized as a mass carried
by the interacting Brownian particles are described
\endabstract

\endtopmatter
\rightheadtext{
ONE BROWNIAN STOCHASTIC FLOW
}
\leftheadtext{
Andrey A.Dorogovtsev
}

\document
\head
Introduction
\endhead
Let $W$  be an
$\mbR^d$-valued Wiener sheet on $\mbR^d\times[0;1].$  Assume that
$\vf\in C^\infty_0(\mbR^d)$  is the spherically symmetric nonnegative
function with the property
$$
\int_{\mbR^d}\vf(u)du=1.
$$
Define for $\ve>0$
$$
\vf_\ve(u)=\ve^{-\frac{d}{2}}\vf(\ve^{-1}u)^{\frac{1}{2}}, u\in\mbR^d.
$$
Note that $\vf_\ve\in C^\infty_0(\mbR^d)$  under every $\ve>0.$  Let us
consider now the equation
$$
\cases
dx_\ve(u,t)=\int_{\mbR^d}\vf_\ve(x_\ve(u,t)-q)W(dq,dt),\\
x_\ve(u,0)=u.
\endcases
\tag1
$$
The solution $x_\ve$  has two important properties. The first one
is that $x_\ve$ is the flow of the gomeomorphisms \cite{1} and the
second one is that for every $u\in\mbR^d$  $\{x_\ve(u,t);
t\geq0\}$  is the Wiener process. Really, denote by $\{\cF_t;
t\geq0\}$   the flow of $\sigma$-field generated by $W$  in a
usual way. Then $\{x_\ve(u,t); t\geq0\}$  is a continuous
$\cF_t$-martingale with the matrix characteristics
$$
\langle
x_\ve(u,\cdot)\rangle_t=\int^t_0\int_{\mbR^d}\vf^2_\ve(x_\ve(u,s)-q)dqdsI=tI,
$$
$I$ is the identity matrix. Hence $\{x_\ve(u,t); t\geq0\}$  is the
Wiener process. Note that for different  $u_1,u_2\in\mbR^d$
$x_\ve(u_1,\cdot)$ and $x_\ve(u_2,\cdot)$ are not independent.
Their joint characteristic equals to
$$
\langle x_\ve(u_1,\cdot), x_\ve(u_2,\cdot)\rangle_t=
$$
$$
= \int^t_0\int_{\mbR^d} \vf_\ve(x_\ve(u_1,s)-q)
\vf_\ve(x_\ve(u_2,s)-q)dqdsI.
$$
So the flow $x_\ve$   now consists of the Wiener processes which
do not stick together. Let us mention that the support $\vf_\ve$
tends to the origin when $\ve\to0+.$ So, one can expect that in
the limit $x_\ve$  turns into the family of independent Wiener
processes. But from another side if $d=1$ and $u_1<u_2,$ then, as
it was mentioned above, $x_\ve(u_1,t)<x_\ve(u_2,t)$ for every $t$
with probability one. Consequently, on the non-formal level
$x_\ve$  in the limit turns into the family of the Wiener
particles, which start from every point of the space and move
independently up to the meeting after which continue the motion
together. The formal realization of this idea on the level of the
description of the particle motion meets some technical troubles
(see \cite{2}). So we propose here to speak not about the
particles but about the mass which they carry. In another words we
will consider  the measure-valued processes related to the flow
$x_\ve$ and their weak limit under $\ve\to0+.$   Let us fix the
probability measure $\mu_0.$ Define the measure-valued process
$\{\mu^\ve_t; t\in[0;1]\}$ in the following way
$$
\mu^\ve_t=\mu_0\circ x_\ve(\cdot,t)^{-1}.\tag2
$$
We will consider this process in the different spaces of measures.
In order to introduce this spaces define for $u\in\mbR^d$
$$
\vf_0(u)=\frac{\|u\|}{1+\|u\|},
$$
$$
\vf_n(u)=\|u\|^n, \ n\geq1.
$$
Denote by $\mfM_n$  the space of all probability measures $\mu$ on
$\mbR^d$ which have the property
$$
\int_{\mbR^d}\vf_n(u)\mu(du)<+\infty.
$$
Note that $\mfM_0$ contains all probability measures. Define the
Wassershtain distance on $\mfM_n$   by the rule \newline
$\forall\mu,\nu\in\mfM_n:$
$$
\gamma_n(\mu,\nu)^{n\vee1}=\inf_{\vk\in C(\mu,\nu)}
\iint_{\mbR^d}\vf_n(u-v)\vk(du,dv),
$$
where the infimum is taken over all probability measures $\vk$ on
$\mbR^d\times\mbR^d$ which have $\mu$  and $\nu$  as their marginal
distributions. It is known, that $(\mfM_n,\gamma_n)$ is the complete
separable metric space \cite{3}  for every $n\geq0.$ The convergence in
$\gamma_0$ is equivalent to the weak convergence. Since $x_\ve$ is continuous
with respect to both variables then $\mu^\ve$ is continuous process in
$\mfM_0$ with probability one. Here we will consider the weak
limits of $\mu^\ve$ under $\ve\to0+.$

\head Main result \endhead
In what follows we will need in the
conditions for the weak compactness of the measure-valued
processes from \cite{4}. Here we present the necessary notations
and statements. Let $C_b(\mbR^d)$   and $C_0(\mbR^d)$   be the
spaces of bounded and finite continuous functions on $\mbR^d.$
Consider the sequence $\{f_n; n\geq1\}$   from $C_0(\mbR^d)$ such
that:

1) $\forall n\geq1: \max_{\mbR^d}|f_n|\leq1,$

2) for every $\vf_n\in C_0(\mbR^d)$ which is bounded by 1 there is
the subsequence $\{f_{n_k}; k\geq1\}$  such that
$$
\max_{\mbR^d}|\vf-f_{n_k}|\to0, \ k\to\infty.
$$
Now let $\{g_n; n\geq1\}$ be the sequence from $C_b(\mbR^d)$  such, that:

1) $\forall x\in\mbR^d, n\geq1: 0\leq g_n(x)\leq1,$

2) $\forall x\in\mbR^d, \|x\|\leq n: \ g_n(x)=0,$

3) $\forall x\in\mbR^d, \|x\|\geq n\geq1: \ g_n(x)=1.$

The following statement was proved in \cite{4}.
\proclaim{Theorem 1}
The family $\{\xi_\alpha; \alpha\in\cU\}$  of random elements in
$C([0;1],\mfM_0)$  is weakly relatively compact if and only if the
following conditions hold:

1) for every $k\geq1$   the set of the random processes
$\{\langle\xi_\alpha, f_k\rangle; \alpha\in\cU\}$ is weakly
relatively compact in $C([0;1]),$

2) the set $\{\langle\xi_\alpha, g_k\rangle; \alpha\in\cU,
k\geq1\}$ is weakly relatively compact in $C([0;1]),$

3) $\forall t\in[0;1] \forall\ve>0:$
$$
\sup_{\alpha\in\cU}P\{\langle\xi_\alpha(t),g_k\rangle>\ve\}\to0,
k\to\infty.
$$
\endproclaim

Applying this theorem to the our processes $\{\mu^\ve\}$ we can
get the following statement. \proclaim {Theorem 2} Let the initial
measure $\mu_0\in\mfM_n,$ $n>2.$  Then the family $\{\mu^\ve;
\ve>0\}$  is weakly compact in $C([0;1],\mfM).$
\endproclaim
\demo{Proof} Let us use the criterion of the weak compactness for
the measure-valued processes from theorem 1. Take the functions
$\{f_k; k\geq1\}$ and $\{g_k; k\geq1\}$ in such a way, that

1) for every $k\geq1, f_k$   and $g_k$  satisfy the Lipshitz condition,

2) the Lipshitz constant for $g_k$ equals 1 for every $k\geq1.$

Now check the weak compactness of the families $\{\langle f_k,
\mu^\ve\rangle; \ve>0\}$ and $\{\langle g_k, \mu^\ve\rangle;
\ve>0, k\geq1\}$ in $C([0;1]).$  Consider the function $h$  on
$\mbR^d$ which satisfies the Lipshitz condition with the constant
$C.$  For the such function
$$
E|\langle h,\mu^\ve_{t_1}\rangle-\langle
h,\mu^\ve_{t_2}\rangle|^n=
$$
$$
=E\bigg|\int_{\mbR^d}
(h(x_\ve(u,t_1))-h(x_\ve(u,t_2)))\mu_0(du)
\bigg|^n\leq
$$
$$
\leq C^n\int_{\mbR^d}
E\|x_\ve(u,t_1)-x_\ve(u,t_2)\|^n\mu_0(du)\leq
C^n\cdot K_n\cdot|t_2-t_1|^{\frac{n}{2}},
$$
where the constant $K_n$  depends on only from $n$  and the dimension
$d.$   This estimation together with relation
$$
\lim_{k\to\infty}\langle g_k,\mu_0\rangle=0
$$
gives us that the conditions 1) and 2) of the theorem 1 hold. In
order to check the condition 3) consider
$$
E
\int_{\mbR^d}\|u\|^n\mu^\ve_t(du)=
E\int_{\mbR^d}\|u+(x_\ve(u,t)-u)\|^n\mu^\ve_t(du)\leq
$$
$$
\leq D
\bigg(
\int_{\mbR^d}\|u\|^n\mu_0(du)+1\bigg),
$$
where the constant $D$  depends only on the dimension $d$  and $n.$  Then
$$
\varlimsup_{k\to\infty}\sup_{\ve>0} P\{\langle
g_k,\mu^\ve_t\rangle
>\delta\}\leq
$$
$$
\leq\varlimsup_{k\to\infty}\sup_{\ve>0}\frac{1}{\delta}E \langle
g_k,\mu^\ve_t\rangle\leq
$$
$$
\leq\varlimsup_{k\to\infty}\sup_{\ve>0}
\frac{1}{\delta}
\frac{1}{k^n}E
\int_{\mbR^d}\|u\|^n\mu^\ve_t(du)\leq
$$
$$
\leq \lim_{k\to\infty}
\frac{1}{\delta}
\frac{1}{k^n}D
\bigg(
\int_{\mbR^d}\|u\|^n\mu_0(du)+1\bigg)=0
$$
for every $t\in[0;1]$  and $\delta>0.$  So the condition 3) holds
and the theorem is proved.
\enddemo

In order to understand how many limit points the family
$\{\mu^\ve\}$ admits under $\ve\to0+$  let us study the behaviour
of the finite-point processes
$\{\vec{x}_\ve(t)=(x_\ve(u_1,t),\ldots,x_\ve(u_n,t)); t\in[0;1]\}$
under $\ve\to0+.$   From now we will consider the case $d=1.$  We
begin with the construction of the Markov process in $\mbR^n$
which can serve as a weak limit of $\vec{x}_\ve.$ Nonformally this
process can be described as follows. In the space $\mbR^n$ we
consider the usual Wiener process up to the first time when the
some of its coordinates became to be equal. After this time
process turns into the Wiener process on the hyperplane, where
this coordinates remains equal. This procedure goes on until we
get the one-dimensional Wiener process. From that moment our
process coincide with this Wiener process. In order to construct
such a random process rigorously and check that it is the unique
weak limit of $\vec{x}_\ve$  let us prove the next theorem.
\proclaim{Theorem 3} The family $\{\vec{x}_\ve; \ve>0\}$  weakly
converges under $\ve\to0+$ in the space $C([0;1], \mbR^n).$
\endproclaim
\demo
The weak compactness of $\{\vec{x}_\ve; \ve>0\}$ follows from the arguments
which were mentioned in the proof of the theorem 2.  So
we only have to prove that under $\ve\to0+$ there is only one limit point.
Fix $\ve>0$  and consider for $x\in\mbR$  the function
$$
g_\ve(x)=\int_\mbR\vf_\ve(x+q)\vf_\ve(q)dq.
$$
The process $\vec{x}_\ve$ is a diffusion process in $\mbR^n$  with zero
drift and the diffusion matrix
$$
A(\vec{x})=\big(g_\ve(x_i-x_j)\big)^n_{ij=1},
$$
where $\vec{x}=(x_1,\ldots,x_n).$  It results from the condition on
$\vf_\ve$  that $A$  coincide  with the identity matrix on the set
$$
G_\ve=\{\vec{x}: |x_i-x_j|>2\ve, i\ne j\}.
$$
Let us define the random moment $\tau_\ve$  as the first exit time
from $G_\ve.$ In our case $\vec{x}(0)\in G_\ve$  because we take
the different initial values $u_1,\ldots,u_n.$  It results from
theorem 1.13.2 \cite{5} that the distribution of the process
$\{\vec{x}_\ve(\tau_\ve\wedge t); t\in[0;1]\}$ coincide with the
distribution of the process $\{\vec{w}(\tau_\ve\wedge t);
t\in[0;1]\},$ where $\vec{w}$  is the Wiener process in $\mbR^n$
starting from $(u_1,\ldots,u_n)$  and $\tau_\ve$  has the same
meaning for $\vec{w}$  as for $\vec{x}_\ve.$   Now suppose that
$u_1<\ldots<u_n$  with out loss of generality. Consider in
$C([0;1],\mbR^n)$  the set
$$
\cG_\delta=\{\vec{f}: f_i(0)=u_i, i=1,\ldots,n, \vec{f}(t)\in
G_\delta, t\in[0;1]\}.
$$
Note that $\cG_\delta$   is an open set in $C([0;1],\mbR^n)$ for
sufficiently small $\delta>0$  (we use the usual uniform norm in
the space of continuous functions). Let $\vk$  be the limit point
of the distributions $\vec{x}_\ve$ under $\ve\to0+.$ It results
from the previous considerations that for all sufficiently small
$\ve>0$  the restriction on $\cG_\delta$  of the distribution of
$\vec{x}_\ve$  coincide with the restriction of the Wiener measure
related to the initial value $(u_1,\ldots,u_n).$   Consequently,
$\vk$  coincide with the Wiener measure on the set $\cG_\delta$
for arbitrary $\delta>0.$   Denote by $\cG$  the closure of the
union
$$
\bigcup_{\delta>0}\cG_\delta.
$$
Note that for all $\ve>0$

$$
P\{\vec{x}_\ve\in\cG\}=1.
$$
So, by characterization of the weak convergence \cite{6},
$$
\vk(\cG)=1.
$$
It remains to describe $\vk$  on the boundary of $\cG.$ In order
to do this let us recall that the random process, which is
obtained from $\vec{x}_\ve$  by choosing some of its coordinates
has the same properties as $\vec{x}_\ve$  under $\ve\to0+.$  Hence
the measure $\vk$  has the following properties. The every
coordinate has the Wiener distribution. Any two coordinates move
as an independent Wiener processes up to their meeting and move
together after this moment. Now the uniqueness of $\vk$   can be
obtained by induction. Theorem is proved.
\enddemo

This theorem has the following consequence. \proclaim{Corollary}
The measure-valued processes $\{\mu^\ve\}$   constructed in
theorem 2   converge weakly under $\ve\to0+.$
\endproclaim
\demo{Proof} In view of the theorem 3   we have only to check the
uniqueness of the limit point for $\{\mu^\ve\}$  under $\ve\to0+.$
Let $\{\nu_t; t\in[0;1]\}$   be the measure-valued process
representing the limit point of $\{\mu^\ve\}$  under $\ve\to0+.$
Take $t_1,\ldots,t_d\in[0;1]$  and for bounded continuous
functions $\vf_1,\ldots,\vf_d$   consider the value
$$
E\prod^d_{k=1}\int_\mbR\vf_k(u)\nu_{t_k}(du).
\tag4
$$
Note, that the set of all such values uniquely define the distribution of
$\{\nu_t; t\in[0;1]\}.$

Due to the previous theorem and the Lebesgue dominated convergence theorem
$$
E\prod^d_{k=1}\int_\mbR\vf_k(u)\nu_{t_k}(du)= \lim_{n\to0+}
E\prod^d_{k=1}\int_\mbR\vf_k(u)\mu_{t_k}^{\ve_n}(du)=
$$
$$
= \lim_{n\to0+} E\prod^d_{k=1}\int_\mbR\vf_k( x_{\ve_n}(u,
t_k))\mu_0(du)=
$$
$$
= \lim_{n\to0+} E \mathop{\int\ldots\int}\limits_{\mbR^d}
\prod^d_{k=1} \vf_k(x_{\ve_n}(u_k,
t_k))\mu_0(du_1)\ldots\mu_0(du_k)=
$$
$$
= \mathop{\int\ldots\int}\limits_{\mbR^d} \lim_{n\to0+} E
\prod^d_{k=1} \vf_k(x_{\ve_n}(u_k,
t_k))\mu_0(du_1)\ldots\mu_0(du_k).
$$
Since the limit in the last integral does not depend on the choice of
the sequence $\ve_n\to0+, n\to\infty,$   then the value of (4)
  is uniquely defined. Hence we have only one limit
point for $\{\mu^\ve\}$  under $\ve\to0+$   and the our statement is
proved.
\enddemo

\enddocument